\documentclass[12pt]{article}%
\usepackage{amsmath}
\usepackage{amsfonts}
\usepackage{amssymb}
\usepackage{graphicx}
\usepackage{ifpdf}
\usepackage{graphicx,amssymb,lineno, epsfig, graphics}%
\providecommand{\U}[1]{\protect\rule{.1in}{.1in}}
\newtheorem{theorem}{Theorem}
\newtheorem{acknowledgement}[theorem]{Acknowledgement}
\newtheorem{remark}[theorem]{Remark}

\begin{document}

\title{Asymptotic analysis of the Bell polynomials by the ray method}
\author{Diego Dominici \thanks{dominici@math.tu-berlin.de}
\\Technische Universit\"{a}t Berlin\\Sekretariat MA 4-5\\Stra\ss e des 17. Juni 136 \\D-10623 Berlin \\Germany
\\Permanent address: \\Department of Mathematics State University of New York at New Paltz \\1 Hawk Dr. \\New Paltz, NY 12561-2443 \\USA}

\maketitle

\begin{abstract}
We analyze the Bell polynomials $B_{n}(x)$  asymptotically as $n\rightarrow\infty$. We obtain asymptotic approximations from the
differential-difference equation which they satisfy, using a discrete version of the ray method. We give some examples showing the accuracy of our formulas.
\end{abstract}
{Keywords Bell polynomials, asymptotic expansions, Stirling numbers}
MSC-class: 34E05, 11B73, 34E20

\section{Introduction}

The Bell polynomials $B_{n}(x)$ are defined by \cite{MR1503161}%
\[
B_{n}(x)=\sum_{k=0}^{n}S_{k}^{n}x^{k},\quad n=0,1,\ldots,
\]
where $S_{k}^{n}$ is a Stirling number of the second kind \cite[24,1,4]%
{MR1225604}. They have the generating function%
\begin{equation}
\sum_{n=0}^{\infty}B_{n}(x)\frac{t^{n}}{n!}=\exp\left[  x\left(
e^{t}-1\right)  \right]  , \label{gener}%
\end{equation}
from which it follows that
\begin{equation}
B_{0}(x)=1 \label{initial}%
\end{equation}
and%
\begin{equation}
B_{n+1}(x)=x\left[  B_{n}^{\prime}(x)+B_{n}(x)\right]  ,\quad n=0,1,\ldots.
\label{diffdiff}%
\end{equation}

The asymptotic behavior of $B_{n}(x)$ was studied by Elbert \cite{MR1820893},
\cite{MR1820892} and Zhao \cite{MR1947752}, using the saddle point method and
(\ref{gener}). In this paper we will use a different approach and analyze
(\ref{diffdiff}) instead of (\ref{gener}). The advantage of our method is that
no knowledge of a generating function is required and therefore it can be
applied to other sequences of polynomials satisfying differential-difference
equations \cite{hermitedif}, \cite{hermitegen}.

\section{Asymptotic analysis}

To analyze (\ref{diffdiff}) asymptotically as $n\rightarrow\infty,$ we use a
discrete version of the ray method \cite{MR1276912}. Replacing the anszat
\begin{equation}
B_{n}(x)=\varepsilon^{-n}F\left(  \varepsilon x,\varepsilon n\right)
\label{BF}%
\end{equation}
in (\ref{diffdiff}), we get%
\begin{equation}
F(u,v+\varepsilon)=u\left(  \varepsilon\frac{\partial F}{\partial x}+F\right)
, \label{f1}%
\end{equation}
with%
\begin{equation}
u=\varepsilon x,\quad v=\varepsilon n \label{u,v}%
\end{equation}
and $\varepsilon$ is a small parameter. We consider asymptotic solutions for
(\ref{f1}) of the form%
\begin{equation}
F(u,v)\sim\exp\left[  \varepsilon^{-1}\psi\left(  u,v\right)  \right]  K(u,v),
\label{f2}%
\end{equation}
as $\varepsilon\rightarrow0.$ Using (\ref{f2}) in (\ref{f1}) we obtain, to
leading order, the eikonal equation
\begin{equation}
e^{q}-u\left(  p+1\right)  =0 \label{eikonal}%
\end{equation}
and the transport equation%
\begin{equation}
\frac{\partial K}{\partial v}+\frac{1}{2}\frac{\partial^{2}\psi}{\partial
v^{2}}K-u\exp\left(  -\frac{\partial\psi}{\partial v}\right)  \frac{\partial
K}{\partial u}=0, \label{transport}%
\end{equation}
where
\begin{equation}
p=\frac{\partial\psi}{\partial x},\quad q=\frac{\partial\psi}{\partial v}.
\label{pq}%
\end{equation}
The initial condition (\ref{initial}), implies%
\begin{equation}
\psi\left(  u,0\right)  =0,\quad K(u,0)=1. \label{initial2}%
\end{equation}

To solve (\ref{eikonal}) we use the method of characteristics, which we
briefly review. Given the first order partial differential equation%
\[
\mathfrak{F}\left(  u,v,\psi,p,q\right)  =0,
\]
with $p,q$ defined in (\ref{pq}), we search for a solution $\psi(u,v)$ by
solving the system of \textquotedblleft characteristic
equations\textquotedblright\
\begin{align*}
u  &  =\frac{du}{dt}=\frac{\partial\mathfrak{F}}{\partial p},\quad v=\frac
{dv}{dt}=\frac{\partial\mathfrak{F}}{\partial q},\\
\dot{p}  &  =\frac{dp}{dt}=-\frac{\partial\mathfrak{F}}{\partial u}%
-p\frac{\partial\mathfrak{F}}{\partial\psi},\quad\dot{q}=\frac{dq}{dt}%
=-\frac{\partial\mathfrak{F}}{\partial v}-q\frac{\partial\mathfrak{F}%
}{\partial\psi},\\
\dot{\psi}  &  =\frac{d\psi}{dt}=p\frac{\partial\mathfrak{F}}{\partial
p}+q\frac{\partial\mathfrak{F}}{\partial q},
\end{align*}
where we now consider $\left\{  u,v,\psi,p,q\right\}  $ to all be functions of
the new variables $t$ and $s.$

For (\ref{eikonal}), we have%
\[
\mathfrak{F}\left(  u,v,\psi,p,q\right)  =e^{q}+p-2u
\]
and therefore the characteristic equations are%
\begin{equation}
\dot{u}+u=0,\quad\dot{v}=e^{q},\quad\dot{p}-p=1,\quad\dot{q}=0\label{sysray}%
\end{equation}
Solving (\ref{sysray}), subject to the initial conditions%
\begin{equation}
u(0,s)=s,\quad v(0,s)=0,\quad p(0,s)=B(s)-1,\label{IC}%
\end{equation}
we obtain%
\[
u=se^{-t},\quad v=Bst,\quad p=Be^{t}-1,\quad q=\ln\left(  Bs\right)
\]
where we have used%
\[
0=\left.  \mathfrak{F}\right\vert _{t=0}=e^{q(0,s)}-sB.
\]
From (\ref{initial2}) and (\ref{IC}) we have%
\begin{equation}
\psi(0,s)=0,\quad K(0,s)=1,\label{ft=0}%
\end{equation}
\ which implies%
\begin{align*}
0 &  =\frac{d}{ds}\psi\left(  0,s\right)  =p(0,s)\frac{d}{ds}u\left(
0,s\right)  +q(0,s)\frac{d}{ds}v\left(  0,s\right)  \\
&  =(B-1)\times1+\ln\left(  Bs\right)  \times0=B-1.
\end{align*}
Thus,
\begin{equation}
u=se^{-t},\quad v=st,\quad p=e^{t}-1,\quad q=\ln\left(  s\right)
.\label{rays}%
\end{equation}

The characteristic equation for $\psi$ is
\[
\dot{\psi}=p\dot{u}+q\dot{v}=\left(  e^{t}-1\right)  \left(  -se^{-t}\right)
+\ln\left(  s\right)  s,
\]
which together with (\ref{ft=0}) gives%
\begin{equation}
\psi\left(  t,s\right)  =s\left(  1-t-e^{-t}\right)  +\ln\left(  s\right)  st.
\label{psi}%
\end{equation}

We shall now solve the transport equation (\ref{transport}). From
(\ref{rays}), we get%
\begin{equation}
\frac{\partial t}{\partial u}=-\frac{te^{t}}{s\left(  t+1\right)  },\quad
\frac{\partial t}{\partial v}=\frac{1}{s\left(  t+1\right)  },\quad
\frac{\partial s}{\partial u}=\frac{e^{t}}{t+1},\quad\frac{\partial
s}{\partial v}=\frac{1}{t+1} \label{inverse}%
\end{equation}
and therefore,
\begin{equation}
\frac{\partial^{2}\psi}{\partial v^{2}}=\frac{\partial q}{\partial v}%
=\frac{\partial q}{\partial t}\frac{\partial t}{\partial v}+\frac{\partial
q}{\partial s}\frac{\partial s}{\partial v}=\frac{1}{s\left(  t+1\right)  }.
\label{psivv}%
\end{equation}
Using (\ref{inverse})-(\ref{psivv}) to rewrite (\ref{transport}) in terms of
$t$ and $s$, we have%
\[
\dot{K}+\frac{1}{2\left(  t+1\right)  }K=0
\]
with solution%
\begin{equation}
K(t,s)=\frac{1}{\sqrt{t+1}}, \label{K}%
\end{equation}
where we have used (\ref{ft=0}).

Solving for $t,s$ in (\ref{rays}), we obtain%
\begin{equation}
t=\mathrm{LW}\left(  \frac{v}{u}\right)  ,\quad s=\frac{v}{\mathrm{LW}\left(
\frac{v}{u}\right)  }, \label{t,s}%
\end{equation}
where $\mathrm{LW}\left(  \cdot\right)  $ denotes the Lambert-W function
\cite{MR1414285}, defined by%
\[
\mathrm{LW}\left(  z\right)  \exp\left[  \mathrm{LW}\left(  z\right)  \right]
=z.
\]
Replacing (\ref{t,s}) in (\ref{psi}) and (\ref{K}), we get%
\[
\psi\left(  u,v\right)  =\frac{v}{\mathrm{LW}\left(  \frac{v}{u}\right)
}+v\ln\left[  \frac{v}{\mathrm{LW}\left(  \frac{v}{u}\right)  }\right]
-(u+v),
\]%
\[
K(u,v)=\frac{1}{\sqrt{\mathrm{LW}\left(  \frac{v}{u}\right)  +1}}%
\]
and from (\ref{f2}) we find that%
\begin{equation}
F(u,v)\sim\exp\left\{  \frac{v/\varepsilon}{\mathrm{LW}\left(  \frac{v}%
{u}\right)  }+\frac{v}{\varepsilon}\ln\left[  \frac{v}{\mathrm{LW}\left(
\frac{v}{u}\right)  }\right]  -(\frac{u+v}{\varepsilon})\right\}  \frac
{1}{\sqrt{\mathrm{LW}\left(  \frac{v}{u}\right)  +1}}, \label{F3}%
\end{equation}
as $\varepsilon\rightarrow0.$ Using (\ref{u,v}) and (\ref{F3}) in (\ref{BF}),
we conclude that%
\begin{equation}
B_{n}(x)\sim\exp\left\{  \frac{n}{\mathrm{LW}\left(  \frac{n}{x}\right)
}+n\ln\left[  \frac{n}{\mathrm{LW}\left(  \frac{n}{x}\right)  }\right]
-(x+n)\right\}  \frac{1}{\sqrt{\mathrm{LW}\left(  \frac{n}{x}\right)  +1}},
\label{Basymp}%
\end{equation}
as $n\rightarrow\infty.$

\begin{figure}[ptb]
\begin{center}
$%
\begin{array}
[c]{c@{\hspace{0.1in}}c}%
\multicolumn{1}{l}{\mbox{}} & \multicolumn{1}{l}{\mbox{}}\\
\rotatebox{270}{\resizebox{2.4in}{!}{\includegraphics{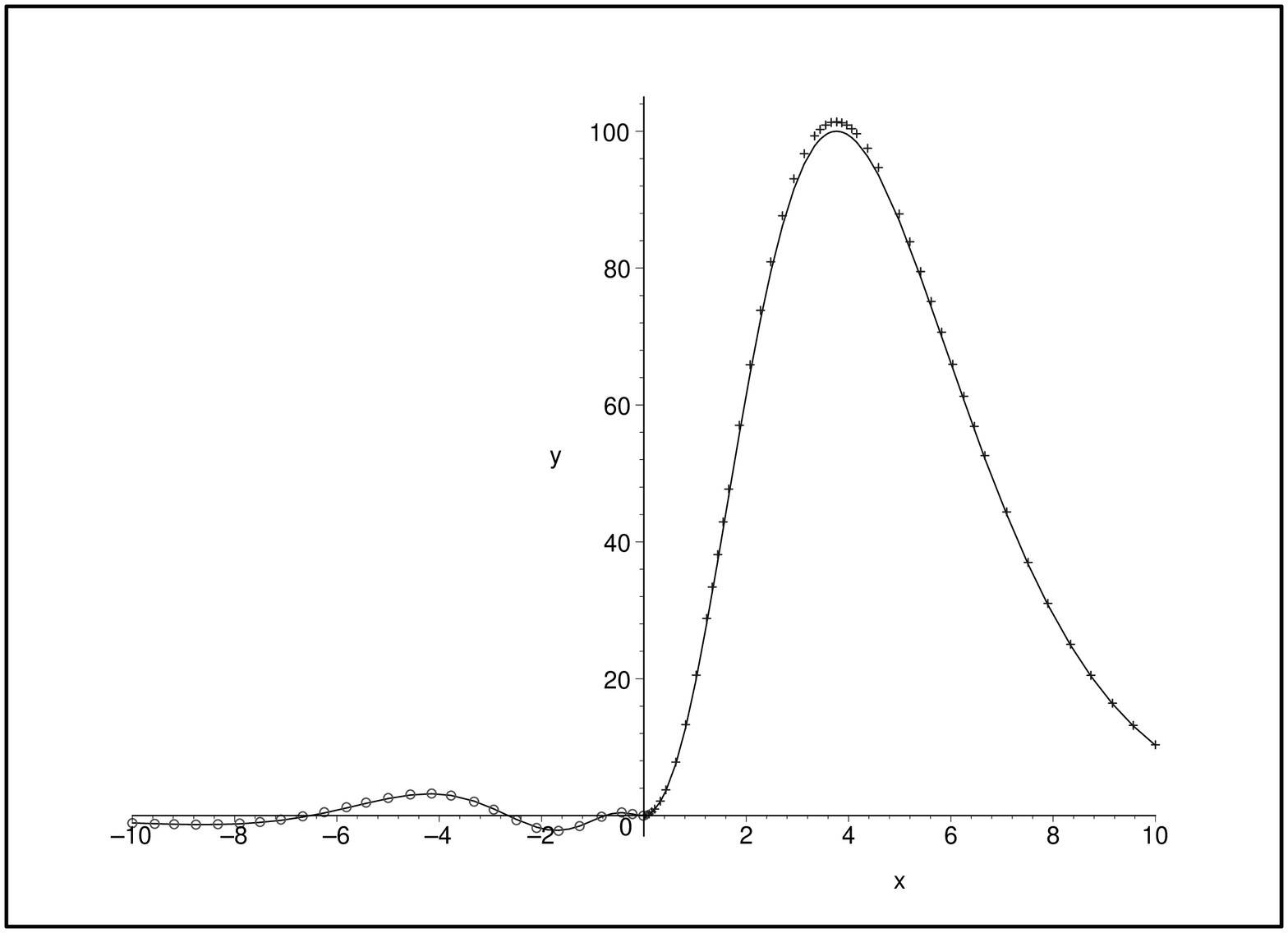}}} &
\rotatebox{270}{\resizebox{2.4in}{!}{\includegraphics{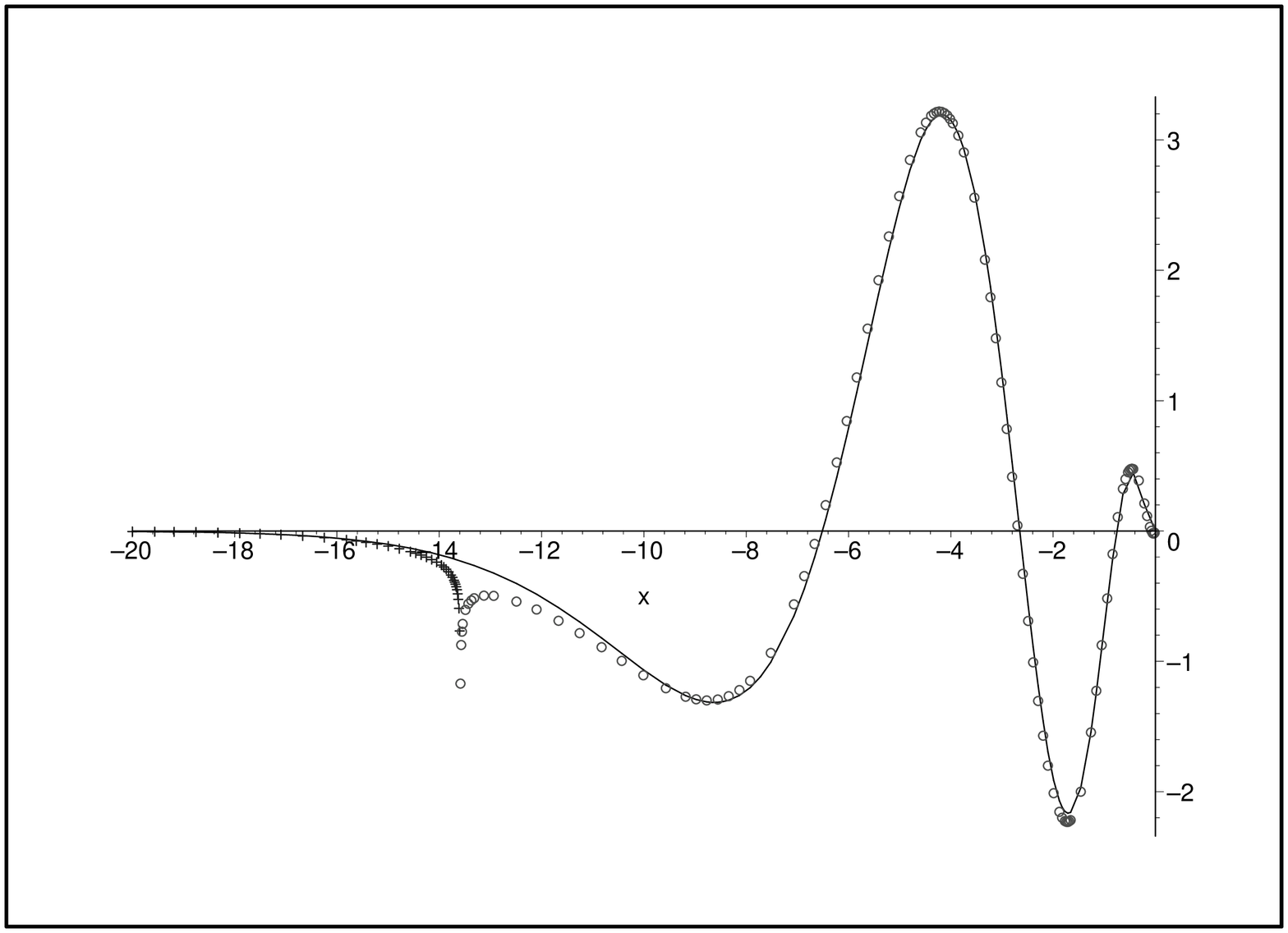}}}\\
& \\
\mbox{(a)} & \mbox{(b)}
\end{array}
$
\end{center}
\caption{A comparison of the exact (solid curve) and asymptotic (ooo), (+++)
values of $B_{5}(x)$.}%
\label{Figure}%
\end{figure}

\begin{remark}
The function $\mathrm{LW}\left(  z\right)  $ has two real-valued branches for
$-e^{-1}\leq z<0,$ denoted by $\mathrm{LW}_{0}\left(  z\right)  $ (the
principal branch of $\mathrm{LW)}$ and $\mathrm{LW}_{-1}\left(  z\right)  ,$
satisfying%
\[
\mathrm{LW}_{0}:\left[  -e^{-1},0\right)  \rightarrow\left[  -1,0\right)
,\quad\mathrm{LW}_{-1}:\left[  -e^{-1},0\right)  \rightarrow\left(
-\infty,-1\right]  ,
\]
with
\[
\mathrm{LW}_{0}\left(  -e^{-1}\right)  =-1=\mathrm{LW}_{-1}\left(
-e^{-1}\right)  .
\]
For $z\geq0,$ $\mathrm{LW}\left(  z\right)  $ has only one real-valued branch%
\[
\mathrm{LW}_{0}:\left[  0,\infty\right)  \rightarrow\left[  0,\infty\right)
\]
and for $z<-e^{-1},$ $\mathrm{LW}_{0}\left(  z\right)  $ and $\mathrm{LW}%
_{-1}\left(  z\right)  $ are complex conjugates. Therefore, for (\ref{Basymp})
to be well defined, we need to consider three separate regions:

\begin{enumerate}
\item An exponential region for $x>0$ or $x<-en.$ Here we have%
\begin{equation}
B_{n}(x)\sim\Phi_{n}(x;0),\quad n\rightarrow\infty, \label{exp}%
\end{equation}
where%
\[
\Phi_{n}(x;k)=\exp\left\{  \frac{n}{\mathrm{LW}_{k}\left(  \frac{n}{x}\right)
}+n\ln\left[  \frac{n}{\mathrm{LW}_{k}\left(  \frac{n}{x}\right)  }\right]
-(x+n)\right\}  \frac{1}{\sqrt{\mathrm{LW}_{k}\left(  \frac{n}{x}\right)  +1}%
}.
\]

\item An oscillatory region for $-en<x<0.$ In this interval,%
\begin{equation}
B_{n}(x)\sim\Phi_{n}(x;0)+\Phi_{n}(x;-1),\quad n\rightarrow\infty.
\label{osci}%
\end{equation}
In Figure \ref{Figure} (a) we plot $B_{5}(x)$ and the asymptotic
approximations (\ref{exp}) (+++) and (\ref{osci}) (ooo), all multiplied by
$e^{-\left\vert x\right\vert }$ for scaling purposes, in the interval $\left(
-10,10\right)  .$ We see that our formulas are quite accurate even for small
values of $n$ and that the transition between (\ref{exp}) and (\ref{osci}) is smooth.

\item A transition region for $x\simeq-en.$ We will analyze this region in the
next section.
\end{enumerate}
\end{remark}

In Figure \ref{Figure} (b) we plot $B_{5}(x)$ and (\ref{exp}) (+++) and
(\ref{osci}) (ooo), all multiplied by $e^{x }$, in the interval $\left(
-20,0\right)  .$ We observe that the approximations (\ref{exp}) and
(\ref{osci}) break down in the neighborhood of $-e5\simeq-13,59$.

\subsection{The transition region}

When $x=-en,$ the quantity $\mathrm{LW}\left(  \frac{n}{x}\right)  +1$
vanishes and (\ref{exp}) is no longer valid. To find an asymptotic
approximation in a neighborhood of $-en,$ we introduce the stretched variable
$\beta$ defined by
\begin{equation}
x=-en-\beta n^{\frac{1}{3}},\quad\beta=O(1). \label{beta}%
\end{equation}
For values of $z$ close to $z_{0}=-e^{-1},$ the Lambert-W function can be
approximated by \cite[(4.22)]{MR1414285}%
\begin{equation}
\mathrm{LW}\left(  z\right)  \sim-1+\sqrt{2e\left(  z-z_{0}\right)  }-\frac
{2}{3}e\left(  z-z_{0}\right)  +\frac{11}{36}\sqrt{2e^{3}\left(
z-z_{0}\right)  ^{3}},\quad z\rightarrow-e^{-1}. \label{LW1}%
\end{equation}
Using(\ref{beta}) in (\ref{LW1}), we have,
\begin{equation}
\mathrm{LW}\left(  \frac{n}{-en-\beta n^{\frac{1}{3}}}\right)  \sim
-1+\sqrt{2e^{-1}\beta}n^{-\frac{1}{3}}-\frac{2}{3}e^{-1}\beta n^{-\frac{2}{3}%
}-\frac{7}{36}\sqrt{2e^{-3}\beta^{3}}n^{-1},\quad\beta\rightarrow0.
\label{LW2}%
\end{equation}
Hence,%
\[
\exp\left\{  \frac{n}{\mathrm{LW}_{k}\left(  \frac{n}{x}\right)  }+n\ln\left[
\frac{n}{\mathrm{LW}_{k}\left(  \frac{n}{x}\right)  }\right]  -(x+n)\right\}
\sim\varphi\left(  \beta,n\right)  ,\quad\beta\rightarrow0,
\]
for $k=0,1$ with $x=-en-\beta n^{\frac{1}{3}}$ and
\begin{equation}
\varphi\left(  \beta,n\right)  =\left(  -1\right)  ^{n}\exp\left\{  \left[
\ln(n)+e-2\right]  n-\left(  e^{-1}-1\right)  \beta n^{\frac{1}{3}}\right\}  .
\label{phi}%
\end{equation}

We now consider solutions for (\ref{diffdiff}) of the form%
\begin{equation}
B_{n}(x)=\varphi\left(  \beta,n\right)  \Lambda\left(  \beta\right)
=\varphi\left[  -\left(  e+\frac{x}{n}\right)  n^{\frac{2}{3}},n\right]
\Lambda\left[  -\left(  e+\frac{x}{n}\right)  n^{\frac{2}{3}}\right]  ,
\label{B1}%
\end{equation}
for some function $\Lambda\left(  \beta\right)  .$ Replacing (\ref{B1}) in
(\ref{diffdiff}) and using (\ref{beta}) we obtain, to leading order%
\[
\Lambda^{\prime\prime}-2e^{-3}\beta\Lambda=0,
\]
with solution%
\begin{equation}
\Lambda\left(  \beta\right)  =C_{1}\mathrm{Ai}\left(  2^{\frac{1}{3}}%
e^{-1}\beta\right)  +C_{2}\mathrm{Bi}\left(  2^{\frac{1}{3}}e^{-1}%
\beta\right)  , \label{lambda}%
\end{equation}
where $\mathrm{Ai}\left(  \cdot\right)  ,$ $\mathrm{Bi}\left(  \cdot\right)  $
are the Airy functions.

To determine the constants $C_{1},C_{2}$ in (\ref{lambda}), we shall match
(\ref{exp}) with (\ref{B1}). Using (\ref{beta}) and (\ref{LW2}) in
(\ref{exp}), we have%
\begin{equation}
B_{n}(x)\sim\varphi\left(  \beta,n\right)  \exp\left(  -\frac{2}{3}\sqrt
{2}e^{-\frac{3}{2}}\beta^{\frac{3}{2}}\right)  \left(  2e^{-1}\beta\right)
^{\frac{1}{4}}n^{-\frac{1}{6}},\quad\beta\rightarrow0^{+}. \label{B3}%
\end{equation}
On the other hand, the Airy functions have the well known asymptotic
approximations \cite[(10.4.59, 10.4.63)]{MR1225604}%
\begin{align*}
\operatorname{Ai}\left(  z\right)   &  \sim\frac{1}{2\sqrt{\pi}}\exp\left(
-\frac{2}{3}z^{\frac{3}{2}}\right)  z^{-\frac{1}{4}},\quad z\rightarrow
\infty,\\
\operatorname{Bi}\left(  z\right)   &  \sim\frac{1}{\sqrt{\pi}}\exp\left(
\frac{2}{3}z^{\frac{3}{2}}\right)  z^{-\frac{1}{4}},\quad z\rightarrow\infty
\end{align*}
and therefore we conclude that
\begin{equation}
C_{1}=\sqrt{\pi}2^{\frac{5}{6}}n^{\frac{1}{6}},\quad C_{2}=0. \label{c1c2}%
\end{equation}

Replacing (\ref{lambda}) and (\ref{c1c2}) in (\ref{B1}), we find that for
$x\simeq-en,$ we have%
\[
B_{n}\left(  x\right)  \sim\sqrt{\pi}2^{\frac{5}{6}}n^{\frac{1}{6}}%
\varphi\left(  \beta,n\right)  \mathrm{Ai}\left(  2^{\frac{1}{3}}e^{-1}%
\beta\right)  ,\quad\quad n\rightarrow\infty.
\]
This concludes the asymptotic analysis of $B_{n}(x)$ for large $n.$

\begin{acknowledgement}
This work was completed while visiting Technische Universit\"{a}t Berlin and
supported in part by a Sofja Kovalevskaja Award from the Humboldt Foundation,
provided by Professor Olga Holtz. We wish to thank Olga for her generous
sponsorship and our colleagues at TU Berlin for their continuous help.
\end{acknowledgement}


\begin{thebibliography}{1}
\expandafter\ifx\csname url\endcsname\relax
  \def\url#1{\texttt{#1}}\fi
\expandafter\ifx\csname urlprefix\endcsname\relax\def\urlprefix{URL }\fi

\bibitem{MR1503161}
E.~T. Bell, Exponential polynomials, Ann. of Math. (2) 35~(2) (1934) 258--277.

\bibitem{MR1225604}
M.~Abramowitz, I.~A. Stegun (Eds.), Handbook of mathematical functions with
  formulas, graphs, and mathematical tables, Dover Publications Inc., New York,
  1992.

\bibitem{MR1820893}
C.~Elbert, Weak asymptotics for the generating polynomials of the {S}tirling
  numbers of the second kind, J. Approx. Theory 109~(2) (2001) 218--228.

\bibitem{MR1820892}
C.~Elbert, Strong asymptotics of the generating polynomials of the {S}tirling
  numbers of the second kind, J. Approx. Theory 109~(2) (2001) 198--217.

\bibitem{MR1947752}
Y.-Q. Zhao, A uniform asymptotic expansion of the single variable {B}ell
  polynomials, J. Comput. Appl. Math. 150~(2) (2003) 329--355.

\bibitem{hermitedif}
D.~E. Dominici, Asymptotic analysis of the {H}ermite polynomials from their
  differential-difference equation  (2007) To appear in the Journal of
  Difference Equations and Applications.

\bibitem{hermitegen}
D.~E. Dominici, Asymptotic analysis of the asymptotic analysis of generalized
  {H}ermite polynomials  (2007) Submitted.

\bibitem{MR1276912}
E.~Giladi, J.~B. Keller, Eulerian number asymptotics, Proc. Roy. Soc. London
  Ser. A 445~(1924) (1994) 291--303.

\bibitem{MR1414285}
R.~M. Corless, G.~H. Gonnet, D.~E.~G. Hare, D.~J. Jeffrey, D.~E. Knuth, On the
  {L}ambert {$W$} function, Adv. Comput. Math. 5~(4) (1996) 329--359.

\end{thebibliography}
\end{document}